# WEIL-ÉTALE MOTIVIC COHOMOLOGY

THOMAS GEISSER*

## 1. Introduction

The purpose of this paper is to study cohomology groups for a new Grothendieck topology, called Weil-étale topology, introduced by Lichtenbaum in [14]. For a variety $X$ over a finite field $\mathbb{F}_q$, an étale sheaf on $X$ corresponds to a sheaf on $\bar X = X \times_{\mathbb{F}_q} \bar{\mathbb{F}}_q$ together with a continuous action of the Galois group $\hat G = \hat{\mathbb{Z}}$. On the other hand, a sheaf for the Weil-étale topology corresponds to an étale sheaf on $\bar X$, together with an action of the group $G \subseteq \hat G$ generated by the Frobenius endomorphism $\varphi$. If we denote the category of Weil-étale sheaves by $\mathcal{T}_G$ and the category of étale sheaves by $\mathcal{T}_{\hat G}$, then there is a morphism of topoi $\gamma : \mathcal{T}_G \to \mathcal{T}_{\hat G}$ from the Weil-étale site to the étale site. The functor $\gamma^*$ is the forgetful functor, and for $U$ étale over $\bar X$ we have $\gamma_*\mathcal{F}(U) = \operatorname{colim}_{H\subseteq H_U} \mathcal{F}(U)^H$, where $H$ runs through the subgroups of finite index of the Galois group contained in the stabilizer $H_U$ of $U$. The first half of the paper is devoted to calculating the derived functor $R\gamma_*\mathcal{F}$. The first main result is

**Theorem 1.1.** *The complex $R\gamma_*\mathcal{F}$ is quasi-isomorphic to the complex of sheaves of continuous $\hat G$-modules whose sections at an étale $U$ over $\bar X$ are given by*

$$\operatorname*{colim}_m \bigoplus_{i=0}^{m-1} \mathcal{F}(\varphi^i U) \xrightarrow{t-1} \operatorname*{colim}_m \bigoplus_{i=0}^{m-1} \mathcal{F}(\varphi^i U).$$

Here the integers in the index set are ordered by divisibility. If we denote the entry in position $\mathcal{F}(\varphi^i U)$ by $f^{(i)}$, then the map in the inverse system is

$$\delta_m^n : \bigoplus_{i=0}^{m-1} \mathcal{F}(\varphi^i U) \to \bigoplus_{i=0}^{nm-1} \mathcal{F}(\varphi^i U), \quad f^{(i)} \mapsto \sum_{j=0}^{n-1} f^{(i+mj)},$$

and $t$ sends $f^{(i)}$ to $(\varphi f)^{(i+1)}$. As a consequence, we get

$$\gamma_*\mathcal{F}(U) \cong \operatorname*{colim}_m \mathcal{F}(U)^{m\mathbb{Z}}$$

$$R^1\gamma_*\mathcal{F}(U) \cong \operatorname*{colim}_{m, N_m^n} \mathcal{F}(U)_{m\mathbb{Z}},$$

where $N_m^n$ is the norm map $f \mapsto \sum_{j=0}^{n-1} \varphi^{-jm} f$, and $R^i\gamma_*\mathcal{F} = 0$ for $i > 1$. In the special case where $\mathcal{F}^\cdot$ is the pull-back $\gamma^*\mathcal{G}^\cdot$ of a bounded above complex of étale sheaves, we get a simpler expression

**Theorem 1.2.** *Let $\mathcal{G}^\cdot$ be a bounded above complex of sheaves in $\mathcal{T}_{\hat G}$. Then there is a quasi-isomorphism*

$$R\gamma_*\mathbb{Z} \otimes^L \mathcal{G}^\cdot \xrightarrow[\sim]{\tau} R\gamma_*(\gamma^*\mathcal{G}^\cdot).$$

* Supported in part by JSPS, the Alfred P. Sloan Foundation, and NSF Grant No. 0070850.

Since $R\gamma_*\mathbb{Z}$ sits in a distinguished triangle $\mathbb{Z} \to R\gamma_*\mathbb{Z} \to \mathbb{Q}[-1]$, this implies that there is a distinguished triangle

$$\mathcal{G}^{\cdot} \to R\gamma_*\gamma^*\mathcal{G}^{\cdot} \to \mathcal{G}^{\cdot} \otimes \mathbb{Q}[-1] \to \mathcal{G}^{\cdot}[1],$$

and the boundary map is induced by the composition $\mathbb{Q}[-1] \to \mathbb{Q}/\mathbb{Z}[-1] \xrightarrow{e} \mathbb{Q}/\mathbb{Z}[0] \xrightarrow{\beta} \mathbb{Z}[1]$, where $e$ is the generator of $\mathrm{Ext}^1_{\hat{G}}(\mathbb{Q}/\mathbb{Z}, \mathbb{Q}/\mathbb{Z})$ and $\beta$ the Bockstein homomorphism. In particular, for a complex of torsion sheaves we have $\mathcal{G}^{\cdot} \cong R\gamma_*\gamma^*\mathcal{G}^{\cdot}$, and for a complex of $\mathbb{Q}$-vector spaces we get $R\gamma_*\gamma^*\mathcal{G}^{\cdot} \cong \mathcal{G}^{\cdot} \oplus \mathcal{G}^{\cdot}[-1]$.

In the second half of the paper we study the hypercohomology $H^i_W(X, \mathbb{Z}(n))$ of the motivic complex for the Weil-étale topology, for a smooth variety $X$ of dimension $d$ over the field $\mathbb{F}_q$. The theorem gives a long exact sequence

$$\ldots \to H^i_{\text{ét}}(X, \mathbb{Z}(n)) \to H^i_W(X, \mathbb{Z}(n)) \to H^{i-1}_{\text{ét}}(X, \mathbb{Q}(n)) \xrightarrow{\delta} H^{i+1}_{\text{ét}}(X, \mathbb{Z}(n)) \to \ldots,$$

where the map $\delta$ is the composition

$$H^{i-1}_{\text{ét}}(X, \mathbb{Q}(n)) \to H^{i-1}_{\text{ét}}(X, \mathbb{Q}/\mathbb{Z}(n)) \xrightarrow{\cdot e} H^i_{\text{ét}}(X, \mathbb{Q}/\mathbb{Z}(n)) \xrightarrow{\beta} H^{i+1}_{\text{ét}}(X, \mathbb{Z}(n)).$$

Rationally, $H^i_W(X, \mathbb{Q}(n)) \cong H^i(X, \mathbb{Q}(n)) \oplus H^{i-1}(X, \mathbb{Q}(n))$, and the cup product with $e \in \mathrm{Ext}^1(\mathbb{Z}, \mathbb{Z})$ is given by multiplication with the matrix $\begin{pmatrix} 0 & 0 \\ 1 & 0 \end{pmatrix}$.

For the cohomological dimension, we show that $H^i_W(X, \mathbb{Z}(n)) = 0$ for $i > 2d+1$ if $n \leq d$, and for $i > n+d+1$ if $n \geq d$. If $X$ is connected and proper, then there is a surjective map $H^{2d}_W(X, \mathbb{Z}(d)) \xrightarrow{\deg} \mathbb{Z}$, and $H^{2d+1}_W(X, \mathbb{Z}(d)) = \mathbb{Z}$.

Given a smooth and proper variety $X$ of dimension $d$ and an integer $n$, there are two conjectures about Weil-étale motivic cohomology. Conjecture $L(X, n)$ of Lichtenbaum states that $H^i_W(X, \mathbb{Z}(n))$ is finitely generated, and finite for $i \neq 2n, 2n+1$, and conjecture $C(X, n)$ states that for any prime $l$,

$$H^i_W(X, \mathbb{Z}(n)) \otimes \mathbb{Z}_l \cong H^i_{cont}(X, \mathbb{Z}_l(n)).$$

In particular, $C(X, n)$ implies that Weil-étale motivic cohomology is an integral model for $l$-adic (and $p$-adic) cohomology. If $T(X, n)$ is the conjunction of the Tate's conjecture for $X$ in degree $n$ together with Beilinson's conjecture that rational and numerical equivalence on $X$ agree with rational coefficients in codimension $n$, then we use arguments of Kahn [10, 12] to show

**Theorem 1.3.** *Let $X$ be a smooth projective variety over $\mathbb{F}_q$, and $n$ an integer. Then*

$$C(X, n) + C(X, d-n) \Rightarrow L(X, n) \Rightarrow C(X, n) \Rightarrow T(X, n).$$

*Conversely, if $T(X, n)$ holds for* all *smooth and projective varieties over $\mathbb{F}_q$ and all $n$, then $C(X, n)$ and $L(X, n)$ hold for all $X$ and $n$.*

As anticipated by Lichtenbaum, Weil-étale motivic cohomology can be used to give formulas for special values of $\zeta$-functions of varieties over finite fields. For a bounded complex with finite cohomology groups define

$$\chi(C^{\cdot}) := \prod_i |H^i(C^{\cdot})|^{(-1)^i}$$

and let

$$\chi(X, \mathcal{O}_X, n) = \sum_{i \leq n, j \leq d} (-1)^{i+j}(n-i) \dim H^j(X, \Omega^i).$$

Comparing to the formula for $l$-adic cohomology of Milne [15], we get

**Theorem 1.4.** *Let $X$ be a smooth projective variety such that $C(X,n)$ holds. Then the order $\rho_n$ of the pole of $\zeta(X,s)$ at $s = n$ is rank $H^{2n}_W(X, \mathbb{Z}(n))$, and*

$$\zeta(X,s) = \pm(1-q^{n-s})^{-\rho_n} \cdot \chi(H^*_W(X, \mathbb{Z}(n)), e) \cdot q^{\chi(X, \mathcal{O}_X, n)} \quad as\ s \to n. \quad (1)$$

*If furthermore $C(X, d-n)$ holds, then*

$$\chi(H^*_W(X, \mathbb{Z}(n)), e) = \prod_i |H^i_W(X, \mathbb{Z}(n))_{tor}|^{(-1)^i} \cdot R^{-1} \quad (2)$$

*where $R$ is the determinant of the pairing*

$$H^{2n}_W(X, \mathbb{Z}(n)) \times H^{2(d-n)}_W(X, \mathbb{Z}(d-n)) \to H^{2d}_W(X, \mathbb{Z}(d)) \to \mathbb{Z}.$$

To give explicit evidence, we show that the surjectivity of the cycle map $\mathrm{Pic} \otimes \mathbb{Q}_l \to H^2_{cont}(X, \mathbb{Q}_l(1))$ implies $C(X, 1)$, giving part b) of the theorem below. Using the method of Soulé [19], we also show that $C(X, n)$ holds for a smooth projective variety $X$ of dimension $d$, which can be constructed out of products of smooth projective curves by union, base extension and blow-ups, and for $n \leq 1$ or $n \geq d-1$, giving part c) of the theorem below.

**Theorem 1.5.** *i) $C(X, n)$ holds for $n \leq 0$.*

*ii) $C(X, 1)$ holds for Hilbert modular surfaces, Picard modular surfaces, Siegel modular threefolds, and in characteristic at least $5$ for supersingular and elliptic K3 surfaces.*

*iii) For $n \leq 1$ or $n \geq d-1$, $C(X, n)$ holds for products of curves, abelian varieties, unirational varieties of dimension at most $3$, and for Fermat hypersurfaces.*

In [12], Kahn shows that conjecture $C(X, n)$ is true for $X$ the product of elliptic curves and any $n$.

This paper is based on ideas of Lichtenbaum [14] and Kahn [10]. We wish to thank B. Kahn, S. Lichtenbaum and T. Saito for several helpful comments. The paper was written while the author was visiting the University of Tokyo, which provided excellent working conditions.

## 2. Profinite completion

Throughout this paper, we let $\mathbb{F}_q$ be a finite field, $\bar{\mathbb{F}}_q$ be the algebraic closure of $\mathbb{F}_q$, and $\varphi$ the Frobenius endomorphism $x \to x^q$ of $\bar{\mathbb{F}}_q$ over $\mathbb{F}_q$. The Galois group $\hat{G}$ of $\bar{\mathbb{F}}_q/\mathbb{F}_q$ is isomorphic to the profinite completion $\varprojlim \mathbb{Z}/m$ of $\mathbb{Z}$, and we let $G$ be the subgroup of $\hat{G}$ generated by $\varphi$. Of course, $G$ is isomorphic to $\mathbb{Z}$, but we want to avoid confusion of $G$-modules with abelian groups. The fixed field of $mG$ and of of $m\hat{G}$ is $\mathbb{F}_{q^m}$.

Let $X$ be a scheme of finite type over $\mathbb{F}_q$ and $\bar{X} = X \times_{\mathbb{F}_q} \bar{\mathbb{F}} \xrightarrow{\pi} X$. Given a sheaf $\mathcal{F}$ on $\bar{X}$, we say that $\hat{G}$ acts on $\mathcal{F}$ if for each $g \in \hat{G}$ there is an isomorphism $\sigma(g) : \mathcal{F} \to g^*\mathcal{F}$ satisfying $\sigma(gh) = \sigma(g)\sigma(h)$. For $f \in \mathcal{F}(U)$, we will abreviate $\sigma(g)f \in \mathcal{F}(gU)$ by $gf$. If $\hat{G}$ acts on $\mathcal{F}$, then for each $V$ étale over $X$, $\hat{G}$ acts on $\mathcal{F}(V \times_{\mathbb{F}_q} \bar{\mathbb{F}}_q)$. Moreover, since every étale scheme $U/\bar{X}$ can be defined over some $\mathbb{F}_{q^m}$, i.e. $U = U' \times_{\mathbb{F}_{q^m}} \bar{\mathbb{F}}_q$, the stabilizer $H_U$ of $U$ is of finite index in $G$, and acts on $\mathcal{F}(U)$. We say that $\hat{G}$ acts continuously on $\mathcal{F}$ if for each étale $U$ over $\bar{X}$, the stabilizer $H_U$ acts continuously on $\mathcal{F}(U)$, where we consider $\mathcal{F}(U)$ with the discrete topology.

By Deligne, SGA 7 XIII, 1.1.3, there is an equivalence of categories between étale sheaves on $X$ and étale sheaves on $\bar X$ together with a continuous action of $\hat G$. Under this equivalence, a sheaf $\mathcal{F}$ on $\bar X$ corresponds to the sheaf $\pi_*^{\hat G}\mathcal{F}$, sending an étale $V$ over $X$ to $\mathcal{F}(V\times_{\mathbb{F}_q}\bar{\mathbb{F}}_q)^{\hat G}$. On the other hand, a sheaf $\mathcal{G}$ on $X$ corresponds to $\pi^*\mathcal{G}$ on $\bar X$. We will denote the topos of étale sheaves on $\bar X$ equipped with a continuous action of $\hat G$ by $\mathcal{T}_{\hat G}$.

In [14], Lichtenbaum defines the Weil-étale topology on $X$. He shows that a Weil-étale sheaf on $X$ is equivalent to an étale sheaf on $\bar X$ together with a $G$-action, where $n \in G$ acts on $\bar X$ via $\varphi^n$ and on $\mathcal{F}$ via $\sigma(\varphi^n): \mathcal{F}(U) \to \mathcal{F}(\varphi^n U)$. We will denote the category of sheaves on $\bar X$ equipped with an action of $G$ by $\mathcal{T}_G$.

Consider the functor $\gamma_* = \pi^*\pi_*^G$ from $\mathcal{T}_G$ to $\mathcal{T}_{\hat G}$. It makes a sheaf of $G$-modules into a sheaf of continuous $\hat G$-modules. Given an étale $U$ over $\bar X$, let $H_U$ be the stabilizer of $U$. Then one can check that

$$\gamma_*\mathcal{F}(U) = \operatorname*{colim}_{H\subseteq H_U} \mathcal{F}(U)^H,$$

where $H$ runs through the subgroups of finite index in $G$ which are contained in $H_U$. Since those subgroups are cofinal in the set of all subgroups of finite index, we will write by abuse of notation $\gamma_*\mathcal{F} = \operatorname{colim}_H \mathcal{F}^H$, remembering that even though $\mathcal{F}^H$ is not defined, the colimit is. Since the invariant functor is left exact, $\gamma_*\mathcal{F}$ is a sheaf if $\mathcal{F}$ is. By construction, $\gamma_*\mathcal{F}$ has an action of $\hat G$. More explicitly, given $U$ with stabilizer $H_U$, and $H \subseteq H_U$, $g \in \hat G$ acts as $\sigma(g) = \sigma(\varphi)^i : \mathcal{F}(U)^H \to \mathcal{F}(gU)^H$, if $g = \varphi^i \mod H$. It is easy to check that this is compatible with the inclusion $\mathcal{F}(U)^H \hookrightarrow \mathcal{F}(U)^{H'}$ for $H' \subseteq H$.

**Proposition 2.1.** *There is a morphism of topoi $\gamma: \mathcal{T}_G \to \mathcal{T}_{\hat G}$. The functor $\gamma^*$ is the forgetful functor, and $\gamma_* = \pi^*\pi_*^G$ sends a sheaf $\mathcal{F}$ to $\operatorname{colim}_H \mathcal{F}^H$. In particular, $\gamma_*$ is left exact and preserves injectives. The adjoint transformation $\mathrm{id} \to \gamma_*\gamma^*$ is an isomorphism.*

*Proof.* Let $\mathcal{F}$ be a sheaf with $G$-action and $\mathcal{G}$ be a sheaf with continuous $\hat G$-action. Then $\mathcal{G} = \operatorname{colim}_H \mathcal{G}^H$, and the map

$$\operatorname{Hom}_G(\gamma^*\mathcal{G},\mathcal{F}) \to \operatorname{Hom}_{\hat G}(\mathcal{G},\gamma_*\mathcal{F})$$
$$\alpha \to \operatorname*{colim}_H \alpha|_{\mathcal{G}^H}$$

is clearly an isomorphism with inverse "composition with the adjoint inclusion $\gamma^*\gamma_*\mathcal{F} \to \mathcal{F}$". The fact $\mathcal{F} \cong \gamma_*\gamma^*\mathcal{F}$ follows from the explicit description of $\gamma_*$ and $\gamma^*$.  □

Given two sheaves $\mathcal{F}$ and $\mathcal{G}$ with a $G$-action, the sheaf $\mathcal{H}\mathfrak{om}(\mathcal{G},\mathcal{F})$ is equipped with a $G$-action by $f^g = \sigma(g) \circ f \circ \sigma(g)^{-1}$. If we apply the previous discussion to $\mathcal{H}\mathfrak{om}(\mathcal{G},\mathcal{F})$, we get the sheaf

$$U \mapsto \gamma_*\mathcal{H}\mathfrak{om}(\mathcal{G},\mathcal{F})(U) = \operatorname*{colim}_{H\subseteq H_U} \operatorname{Hom}(\mathcal{G}|_U,\mathcal{F}|_U)^H,$$

where the latter are the homomorphisms which are compatible with the action of $H$. In the special case when $\mathcal{G} = A$ is constant, then by adjointness of global section and constant sheaf functor $\operatorname{Hom}(A,\mathcal{F}|_U)^H = \operatorname{Hom}(A,\mathcal{F}(U))^H$, and the formula simplifies to $\gamma_*\mathcal{H}\mathfrak{om}(A,\mathcal{F})(U) = \operatorname{colim}_{H\subseteq H_U} \operatorname{Hom}(A,\mathcal{F}(U))^H$. In particular, $\gamma_*\mathcal{F} = \gamma_*\mathcal{H}\mathfrak{om}(\mathbb{Z},\mathcal{F})$.

Consider $\oplus_{i=0}^{m-1}\mathcal{F}(\varphi^i U)$, if $mG$ is contained in the stabilizer of $U$. Given an integer $i \mod m$, we write $f^{(i)}$ for $f \in \mathcal{F}(\varphi^i U)$ in the $i$th summand.

**Lemma 2.2.** *Let $U$ be étale over $\bar{X}$ and $H = mG \subseteq G$ be a subgroup contained in the stabilizer of $U$.*

*a) For $\mathcal{F}$ in $\mathcal{T}_G$, there are isomorphisms*

$$\operatorname{Hom}(\mathbb{Z}[G], \mathcal{F}(U))^H \xleftarrow[\sim]{\alpha} \mathbb{Z}[G] \otimes_H \mathcal{F}(U) \xrightarrow[\sim]{\beta} \bigoplus_{i=0}^{m-1} \mathcal{F}(\varphi^i U),$$

*where $\beta(\varphi^a \otimes f) = (\varphi^a f)^{(a)} \in \mathcal{F}(\varphi^a U)$ and*

$$\alpha(s \otimes f)(g) = \begin{cases} gsf & \text{if } gs \in H \\ 0 & \text{otherwise.} \end{cases}$$

*b) Under these isomorphisms, the action of $\varphi^j \in G$ via multiplication on $\mathbb{Z}[G]$ on $\operatorname{Hom}(\mathbb{Z}[G], \mathcal{F}(U))^H$ corresponds to the multiplication on $\mathbb{Z}[G]$ on $\mathbb{Z}[G] \otimes_H \mathcal{F}(U)$, and to the map which sends $f^{(i)} \in \mathcal{F}(\varphi^i U)$ to $(\varphi^j f)^{(i+j)} \in \mathcal{F}(\varphi^{i+j} U)$.*

*c) The conjugation action of $g \in \hat{G}$, $\operatorname{Hom}(\mathbb{Z}[G], \mathcal{F}(U))^H \to \operatorname{Hom}(\mathbb{Z}[G], \mathcal{F}(gU))^H$, $F \mapsto g \circ F \circ g^{-1}$, corresponds to the rotation map*

$$g \cdot - : \bigoplus_{i=0}^{m-1} \mathcal{F}(\varphi^i U) \to \bigoplus_{i=0}^{m-1} \mathcal{F}(\varphi^i gU), \quad f^{(i)} \mapsto f^{(i-a)},$$

*if $g \equiv \varphi^a \mod m\hat{G}$.*

*d) Given a second subgroup $H' = mnG \subseteq H = mG$, then the canonical inclusion $\operatorname{Hom}(\mathbb{Z}[G], \mathcal{F}(U))^H \hookrightarrow \operatorname{Hom}(\mathbb{Z}[G], \mathcal{F}(U))^{H'}$ corresponds to the map*

$$\delta_m^n : \bigoplus_{i=0}^{m-1} \mathcal{F}(\varphi^i U) \to \bigoplus_{i=0}^{nm-1} \mathcal{F}(\varphi^i U), \quad f^{(i)} \to \sum_{j=0}^{n-1} f^{(i+mj)}.$$

*This makes sense because $mG$ stabilizes $U$, i.e. $\varphi^i U = \varphi^{i+mj} U$.*

*e) The map $\delta_m^n$ is compatible with the action of $G$ and of $\hat{G}$.*

*Proof.* a), b) are straightforward verifications.

c) First note that if $g \equiv \varphi^a \mod m\hat{G}$, then $gU = \varphi^a U$, so that the map is defined. It is easy to check that the conjugation map $F \mapsto gFg^{-1}$ corresponds under the isomorphism $\alpha$ to the map

$$\mathbb{Z}[G] \otimes_H \mathcal{F}(U) \to \mathbb{Z}[G] \otimes_H \mathcal{F}(gU)$$
$$s \otimes f \mapsto g^{-1}s \otimes gf$$

and this corresponds under the isomorphism $\beta$ to the rotation map.

d) Let $u$ be the inclusion of $H$-invariant maps into $H'$-invariant maps, let $v : \mathbb{Z}[G] \otimes_H \mathcal{F}(U) \to \mathbb{Z}[G] \otimes_{H'} \mathcal{F}(U)$ be the map $s \otimes f \mapsto \sum_{j=0}^{n-1} s\varphi^{jm} \otimes \varphi^{-jm} f$, and let $\delta_m^n$ be the map of the lemma. Then it is easy to check that $u \circ \alpha = \alpha \circ v$ and $\delta_m^n \circ \beta = \beta \circ v$.

e) This is clear because $\alpha$ and $\beta$ are compatible with the action of $G$ and $\hat{G}$. More concretely, for $\varphi^a \in G$ we have

$$\delta_m^n(\varphi^a \cdot f^{(i)}) = \delta_m^n((\varphi^a f)^{(i+a)}) = \sum_{j=0}^{n-1}(\varphi^a f)^{(i+a+mj)} = \sum_{j=0}^{n-1}\varphi^a \cdot f^{(i+mj)} = \varphi^a \cdot \delta_m^n(f^{(i)}).$$

□

Consider the sheaf $\operatorname{colim}_m \bigoplus_{i=0}^{m-1} (\varphi^i)^* \mathcal{F}$, where the index set is ordered by divisibility, and the map in the direct system are the maps $\delta_m^n$. As usual, this is an abbreviation for the sheaf $U \mapsto \operatorname{colim}_{mG \subseteq H_U} \bigoplus_{i=0}^{m-1} \mathcal{F}(\varphi^i U)$. Since the action of $g \in \hat{G}$ is compatible with $\delta_m^n$, we get an action of $\hat{G}$ on $\operatorname{colim}_m \bigoplus_{i=0}^{m-1} \mathcal{F}(\varphi^i U)$.

**Corollary 2.3.** *The functor $\gamma_* \mathcal{H}om(\mathbb{Z}[G], -)^H$ is exact. In particular, its derived functors $R^i \operatorname{colim}_H \mathcal{H}om(\mathbb{Z}[G], -)^H$ are zero for $i > 0$.*

*Proof.* This follows from $\operatorname{colim}_m \mathcal{H}om(\mathbb{Z}[G], \mathcal{F})^{mG} = \operatorname{colim}_H \bigoplus_{i=0}^{m-1} (\varphi^{-i})^* \mathcal{F}$, because the functor $\mathcal{F} \mapsto (\varphi^{-i})^* \mathcal{F}$ and colimits are exact. □

## 3. The functor $R\gamma_*$

**Theorem 3.1.** *The complex $R\gamma_* \mathcal{F}$ is quasi-isomorphic to the complex of sheaves of continuous $\hat{G}$-modules sending an étale $U/\bar{X}$ to*

$$\operatorname{colim}_m \bigoplus_{i=0}^{m-1} \mathcal{F}(\varphi^i U) \xrightarrow{t-1} \operatorname{colim}_m \bigoplus_{i=0}^{m-1} \mathcal{F}(\varphi^i U).$$

*The map $t$ sends $f^{(i)} \in \mathcal{F}(\varphi^i U)$ to $(\varphi f)^{(i+1)} \in \mathcal{F}(\varphi^{i+1} U)$.*

*Proof.* By lemma 2.2 it suffices to show that $R\gamma_* \mathcal{F}$ is quasi-isomorphic to the complex

$$\operatorname{colim}_m \mathcal{H}om(\mathbb{Z}[G], \mathcal{F})^{mG} \xrightarrow{t-1} \operatorname{colim}_m \mathcal{H}om(\mathbb{Z}[G], \mathcal{F})^{mG}.$$

Let $P.$ be the free resolution $0 \to \mathbb{Z}[G] \xrightarrow{t-1} \mathbb{Z}[G] \to 0$ of the constant sheaf $\mathbb{Z}$, and let $\mathcal{F} \to I^\cdot$ be an injective resolution (which is easily seen to exist). We claim that there are quasi-isomorphisms

$$\operatorname{colim}_m \mathcal{H}om(P., \mathcal{F})^{mG} \xrightarrow{\sim} \operatorname{colim}_m \mathcal{H}om(P., I^\cdot)^{mG} \xleftarrow{\sim} \operatorname{colim}_m \mathcal{H}om(\mathbb{Z}, I^\cdot)^{mG}.$$

Indeed, consider the double complex $\operatorname{colim}_m \mathcal{H}om(P_a, I^b)^{mG}$. Taking first horizontal cohomology, we get $\operatorname{colim}_m \mathcal{H}om(\mathbb{Z}, I^b)^{mG}$ for $a = 0$ and zero otherwise by injectivity of $I^b$. Taking now vertical cohomology we get $R^b \gamma_* \mathcal{F}$ for $a = 0$, and zero otherwise. Conversely, if we first take vertical cohomology, we get complexes

$$R^b \big( \operatorname{colim}_m \mathcal{H}om(\mathbb{Z}[\mathbb{Z}], -)^{mG} \big) \mathcal{F} \xrightarrow{t-1} R^b \big( \operatorname{colim}_m \mathcal{H}om(\mathbb{Z}[\mathbb{Z}], -)^{mG} \big) \mathcal{F}$$

concentrated in degree $a = 0, 1$ for each $b$. But by corollary 2.3 the higher derived functors vanish, and the theorem follows. □

To calculate the cohomology explicitly, let $N_m^n : \mathcal{F}(U) \to \mathcal{F}(U)$ be the map $f \mapsto \sum_{j=0}^{n-1} \varphi^{-jm} f$ for $mG \subseteq H_U$. This descends to a map $N_m^n : \mathcal{F}(U)_{m\mathbb{Z}} \to \mathcal{F}(U)_{mn\mathbb{Z}}$, because

$$N_m^n(\varphi^m f) = \sum_{j=0}^{n-1} \varphi^{-(j-1)m} f = N_m^n(f) + \varphi^m f - \varphi^{-(n-1)m} f = N_m^n(f) + (1 - \varphi^{-nm}) \varphi^m f.$$

**Proposition 3.2.** *Let $\mathcal{F}$ be a sheaf in $\mathcal{T}_G$. Then there are isomorphisms*

$$\operatorname*{colim}_{m} \mathcal{F}^{m\mathbb{Z}} \xrightarrow[\sim]{\Delta} \ker t - 1 =: \gamma_* \mathcal{F}$$

$$R^1 \gamma_* \mathcal{F} := \operatorname{coker} t - 1 \xrightarrow[\sim]{S} \operatorname*{colim}_{m, N_m^n} \mathcal{F}_{m\mathbb{Z}}$$

*In particular, $\gamma_* \mathbb{Z} \cong \mathbb{Z}$ and $R^1 \gamma_* \mathbb{Z} \cong \mathbb{Q}$.*

*Proof.* a) Given a sheaf $\mathcal{F}$ in $\mathcal{T}_G$ and $U$ with $mG \subseteq H_U$, we define

$$\Delta_m : \mathcal{F}^{m\mathbb{Z}}(U) \to \bigoplus_{i=0}^{m-1} \mathcal{F}(\varphi^i U), \quad f \mapsto \sum_{i=0}^{m-1} (\varphi^i f)^{(i)}.$$

The map $\Delta_m$ is an isomorphism to the kernel of $t - 1$, because on the one hand

$$(t-1)\Delta_m(f) = \sum_{i=0}^{m-1} (\varphi^{i+1} f)^{(i+1)} - (\varphi^i f)^{(i)} = (\varphi^m f)^{(0)} - f^{(0)} = 0.$$

On the other hand, given $(f_i^{(i)}) \in \bigoplus_{i=0}^{m-1} \mathcal{F}(\varphi^i U)$, $(t-1)(f_i^{(i)}) = 0$ implies that $f_i = \varphi^i f_0$. It is easy to check that $\delta_m^n \circ \Delta_m = \Delta_{mn}$ for $f \in \mathcal{F}^{m\mathbb{Z}}(U)$, hence we get a map

$$\Delta : \operatorname*{colim}_{m} \mathcal{F}^{m\mathbb{Z}}(U) \to \operatorname*{colim}_{m} \bigoplus_{i=0}^{m-1} \mathcal{F}(\varphi^i U).$$

which is an isomorphism to the kernel of $t - 1$.

b) Consider the map

$$S_m : \bigoplus_{i=0}^{m-1} \mathcal{F}(\varphi^i U) \to \mathcal{F}(U) \quad f^{(i)} \mapsto \varphi^{-i} f$$

We have

$$S_{mn}\left(\delta_m^n(f^{(i)})\right) = S_{mn}\left(\sum_{j=0}^{n-1} f^{(i+jm)}\right) = \sum_{j=0}^{n-1} \varphi^{-i-jm} f = N_m^n(\varphi^{-i} f) = N_m^n\left(S_m(f^{(i)})\right),$$

hence we get a surjective map

$$\tilde{S} : \operatorname*{colim}_{m} \bigoplus_{i=0}^{m-1} \mathcal{F}(\varphi^i U) \to \operatorname*{colim}_{m, N_m^n} \mathcal{F}(U).$$

We claim that this map descends to an isomorphism

$$S : \left( \operatorname*{colim}_{m} \bigoplus_{i=0}^{m-1} \mathcal{F}(\varphi^i U) \right) / t - 1 \to \operatorname*{colim}_{m, N_m^n} \mathcal{F}(U)_{m\mathbb{Z}}.$$

First note that

$$S_m\left((t-1)\sum_{i=0}^{m-1} f_i^{(i)}\right) = S_m \sum_{i=0}^{m-1} \left((\varphi f_i)^{(i+1)} - (f_i)^{(i)}\right)$$

$$= S_m\left((\varphi f_{m-1})^{(0)} - (f_{m-1})^{(m-1)}\right) = (1 - \varphi^{-m})(\varphi f_{m-1}) = 0 \in \mathcal{F}_{m\mathbb{Z}}.$$

On the other hand, from the equation $(t-1)\sum_{j=0}^{i-1} (\varphi^{j-i} f)^{(j)} = f^{(i)} - (\varphi^{-i} f)^{(0)}$ we conclude that $S_m(\sum_i f_i^{(i)}) = \sum_i \varphi^{-i} f_i = 0$ implies $\sum_i f_i^{(i)} \sim \sum_i (\varphi^{-i} f_i)^{(0)} = 0$ modulo the image of $t - 1$.

Finally, the map which is multiplication by $\frac{1}{m}$ on the copy of $\mathbb{Q}$ indexed by $m$ induces an isomorphism $\operatorname{colim}_{m, N_m^n} \mathbb{Q} \xrightarrow{\sim} \mathbb{Q}$. □

**Theorem 3.3.** *Let $\mathcal{G}^\cdot$ be a bounded above complex of sheaves in $\mathcal{T}_{\hat{G}}$. Then there is a quasi-isomorphism*
$$R\gamma_*\mathbb{Z} \otimes^L \mathcal{G}^\cdot \xrightarrow[\sim]{\tau} R\gamma_*(\gamma^*\mathcal{G}^\cdot).$$
*In particular, there is a distinguished triangle*
$$\mathcal{G}^\cdot \to R\gamma_*\gamma^*\mathcal{G}^\cdot \to \mathcal{G}^\cdot \otimes \mathbb{Q}[-1].$$

*Proof.* Theorem 3.1 gives us a complex $R\gamma_*\mathbb{Z}$ which consists of *flat* sheaves. Hence we get a quasi-isomorphism of complexes
$$R\gamma_*\mathbb{Z} \otimes^L \mathcal{G}^\cdot(U) \cong \operatorname*{colim}_m \bigoplus_{i=0}^{m-1} \mathcal{G}^\cdot(U) \xrightarrow{t-1} \operatorname*{colim}_m \bigoplus_{i=0}^{m-1} \mathcal{G}^\cdot(U).$$

Here $t$ rotates the summands. On the other hand, we have a quasi-isomorphism
$$R\gamma_*(\gamma^*\mathcal{G}^\cdot)(U) \cong \operatorname*{colim}_m \bigoplus_{i=0}^{m-1} \mathcal{G}^\cdot(\varphi^i U) \xrightarrow{t-1} \operatorname*{colim}_m \bigoplus_{i=0}^{m-1} \mathcal{G}^\cdot(\varphi^i U).$$

Here $t$ acts as in theorem 3.1. We claim that the following map induces an isomorphism of the two complexes:
$$\operatorname*{colim}_m \bigoplus_{i=0}^{m-1} \mathcal{G}^\cdot(U) \hookleftarrow \operatorname*{colim}_m \bigoplus_{i=0}^{m-1} \mathcal{G}^\cdot(U)^{mG} \xrightarrow{\tau} \operatorname*{colim}_m \bigoplus_{i=0}^{m-1} \mathcal{G}^\cdot(\varphi^i U).$$

The inclusion is an isomorphism because each $\mathcal{G}^i$ is a continuous $\hat{G}$-module, and the map $\tau$ sends the entry $f^{(i)} \in \mathcal{G}(U)$ in the middle group to $(\varphi^i f)^{(i)} \in \mathcal{G}^\cdot(\varphi^i U)$ in the right group. Note that we have to work with the intermediate group because $t\tau(f^{(m-1)}) = (\varphi^m f)^{(0)}$, but $\tau t(f^{(m-1)}) = f^{(0)}$. The exact triangle comes from the description of $R\gamma_*\mathbb{Z}$ in proposition 3.2. □

**Corollary 3.4.** *If $\mathcal{G}^\cdot$ is a bounded above complex of torsion sheaves in $\mathcal{T}_{\hat{G}}$, then*
$$\mathcal{G}^\cdot \cong R\gamma_*\gamma^*\mathcal{G}^\cdot.$$
□

If $\mathcal{F}$ is not a $\hat{G}$-module, then $R\gamma_*\mathcal{F}$ cannot be described as easily.
*Example.* If $G$ acts on the constant sheaf $\mathcal{F} = \mathbb{Q}$ as multiplication by $r \neq \pm 1$, then $R\gamma_*\mathcal{F} = 0$. More generally, for any constant sheaf $\mathcal{F}$ of rational vector spaces, $\gamma_*\mathcal{F}$ is the smallest subspace on which $\varphi$ acts as multiplication by some root of unity, and $R^1\gamma_*\mathcal{F}$ is the smallest quotient space where $\varphi$ acts as multiplication by some root of unity.
*Example.* Let $\mathcal{F}$ be the sheaf $\oplus_i A$ for an abelian group $A$, where $\mathbb{Z}$ acts by shifting the factors. Then $\gamma_*\mathcal{F} = 0$, whereas $R^1\gamma_*\mathcal{F} \subseteq \prod_i A$ consists of elements which are periodic for some period. Indeed, for an element of $\mathcal{F}^{m\mathbb{Z}}$, the entries in the sum are $m$-periodic, hence they must be zero. On the other hand, $\mathcal{F}_{m\mathbb{Z}} \cong \oplus_{i=0}^{m-1} A$ and the map $\mathcal{F}_{m\mathbb{Z}} \to \mathcal{F}_{mn\mathbb{Z}}$ is the $n$-fold concatenation map.

## 4. Extensions of $\mathbb{Z}[G]$-modules

We recall some well known results for which we could not find a good reference. Let $e \in \mathrm{Ext}^1_G(\mathbb{Z}, \mathbb{Z}) \cong H^1(\mathbb{Z}, \mathbb{Z}) \cong \mathbb{Z}$ be the class of a generator.

**Lemma 4.1.** *The element $e \in \mathrm{Ext}^1_G(\mathbb{Z}, \mathbb{Z})$ corresponds to the extension $0 \to \mathbb{Z} \to N \to \mathbb{Z} \to 0$, where $N = \mathbb{Z} \oplus \mathbb{Z}$ as an abelian group, and $a \in G$ sends $(r,s)$ to $(r+as, s)$.*

*Proof.* Consider the projective resolution $0 \to \mathbb{Z}[G] \xrightarrow{t-1} \mathbb{Z}[G] \to \mathbb{Z} \to 0$. The generator $s \in \mathrm{Hom}_G(\mathbb{Z}[G], \mathbb{Z}) = \mathbb{Z}$, $s(\sum_i a_i t^i) = \sum_i a_i$, maps under the Bockstein homomorphism to the generator $e$ of $\mathrm{Ext}_G(\mathbb{Z}, \mathbb{Z})$. Hence $e$ corresponds to the push-out

$$\begin{array}{ccccccccc} 0 & \longrightarrow & \mathbb{Z}[G] & \xrightarrow{t-1} & \mathbb{Z}[G] & \longrightarrow & \mathbb{Z} & \longrightarrow & 0 \\ & & {\scriptstyle s}\downarrow & & \downarrow & & \| & & \\ 0 & \longrightarrow & \mathbb{Z} & \longrightarrow & N' & \longrightarrow & \mathbb{Z} & \longrightarrow & 0. \end{array}$$

It suffices to show that there is a short exact sequence

$$0 \to \mathbb{Z}[G] \xrightarrow{(s,1-t)} \mathbb{Z} \oplus \mathbb{Z}[G] \xrightarrow{(u,v)} N = \mathbb{Z} \oplus \mathbb{Z} \to 0,$$

where $u(a, \sum_i a_i t^i) = a + \sum_i a_i i$, and $v(a, \sum_i a_i t^i) = \sum_i a_i$. If $G$ acts on $\mathbb{Z} \oplus \mathbb{Z}$ in the prescribed way, then $(u,v)$ is a $G$-homomorphism, because

$$t \cdot (u,v)(a, \sum_i a_i t^i) = t(a + \sum_i a_i i, \sum_i a_i) = (a + \sum_i a_i(i+1), \sum_i a_i)$$
$$= (u,v)(a, \sum_i a_i t^{i+1}) = (u,v) t \cdot (a, \sum_i a_i t^i).$$

It is easy to see that the sequence is exact on the left and the right, and that

$$u \circ (s, 1-t)(\sum_i a_i t^i) = u(\sum_i a_i, \sum_i (a_i t^i - a_i t^{i+1})) = \sum_i (a_i + a_i i - a_i(i+1)) = 0,$$

and $v \circ (s, 1-t)(\sum_i a_i t^i) = \sum_i (a_i - a_{i-1}) = 0$. Finally, if $x = (a, \sum_i a_i t^i)$ is in the kernel of $(u,v)$, then $\sum_i a_i = 0$ implies that $\sum_i a_i t^i$ is divisible by $1-t$. If $\sum_i a_i t^i = (1-t)\sum_i b_i t^i$, then $u(a, (1-t)\sum_i b_i t^i) = a + \sum_i b_i i - \sum_i b_i(i+1) = 0$ implies that $a = \sum_i b_i$, hence $x = (s, 1-t)(\sum_i b_i t^i)$. □

For every complex of sheaves $C^\cdot$, the element $e \in \mathrm{Ext}^1_G(\mathbb{Z}, \mathbb{Z})$ gives rise to a map in the derived category $e_{C^\cdot} \in \mathrm{Hom}_G(C^\cdot, C^\cdot[1])$, namely the connecting homomorphism of the distinguished triangle

$$C^\cdot \to C^\cdot \otimes N \to C^\cdot \xrightarrow{\beta} C^\cdot[1].$$

This map is natural, i.e. if $\tau: F \to G$ is transformation of triangulated functors, then the following diagram commutes

$$\begin{array}{ccc} F(C^\cdot) & \xrightarrow{F(e_{C^\cdot})} & F(C^\cdot)[1] \\ {\scriptstyle \tau_{C^\cdot}}\downarrow & & {\scriptstyle \tau_{C^\cdot}}\downarrow \\ G(C^\cdot) & \xrightarrow{G(e_{C^\cdot})} & G(C^\cdot)[1] \end{array}$$

In the special case where $C^\cdot$ is a $G$-module $A$, we get a canonical class $e_A \in \mathrm{Ext}^1_G(A, A)$.

**Lemma 4.2.** *The cup product with $e_A$ induces the natural map $H^0(G, A) = A^G \to A_G = H^1(G, A)$ on cohomology.*

*Proof.* The proof of [18, lemma 1.2] works as well in our situation. □

In analogy to the definition of the extension $N$, we get canonical elements in the following sequence of push-outs and pull-backs of extension groups

$$\mathrm{Ext}^1_G(\mathbb{Z}, \mathbb{Z}) \to \mathrm{Ext}^1_G(\mathbb{Z}, \mathbb{Q}) \leftarrow \mathrm{Ext}^1_G(\mathbb{Q}, \mathbb{Q}) \to \mathrm{Ext}^1_G(\mathbb{Q}, \mathbb{Q}/\mathbb{Z}) \leftarrow \mathrm{Ext}^1_G(\mathbb{Q}/\mathbb{Z}, \mathbb{Q}/\mathbb{Z}).$$

The latter two extensions are restrictions of extensions of $\hat{G}$-modules. The image of $e_{\mathbb{Q}/\mathbb{Z}}$ under the canonical map

$$\mathrm{Ext}^1_{\hat{G}}(\mathbb{Q}/\mathbb{Z}, \mathbb{Q}/\mathbb{Z}) \to \mathrm{Ext}^1_{\hat{G}}(\mathbb{Q}, \mathbb{Q}/\mathbb{Z}) \xrightarrow{\beta} \mathrm{Ext}^2_{\hat{G}}(\mathbb{Q}, \mathbb{Z})$$

induced by the canonical projection and the Bockstein homomorphism is calculated by the following pull-back commutative diagram

$$\begin{array}{ccccccccc}
0 & \to & \mathbb{Z} & \to & \mathbb{Q} & \xrightarrow{g} & M & \to & \mathbb{Q} & \to & 0 \\
& & \downarrow & & \downarrow & & \| & & \| & & \\
0 & \to & \mathbb{Q}/\mathbb{Z} & \to & & & M & \to & \mathbb{Q} & \to & 0 \\
& & \| & & & & \downarrow & & \downarrow & & \\
0 & \to & \mathbb{Q}/\mathbb{Z} & \to & & & N \otimes \mathbb{Q}/\mathbb{Z} & \to & \mathbb{Q}/\mathbb{Z} & \to & 0
\end{array}$$

These extensions have been studied by Kahn [10, def. 4.1]. He denotes the complex $\mathbb{Q} \xrightarrow{g} M$ by $\mathbb{Z}^c$, where $\mathbb{Q}$ is in degree 0, $M$ is in degree 1.

**Theorem 4.3.** *There is a quasi-isomorphism of complexes of $\hat{G}$-modules $R\gamma_*\mathbb{Z} \cong \mathbb{Z}^c$. In particular, $R\gamma_*\mathbb{Q} \cong \mathbb{Q} \oplus \mathbb{Q}[-1]$ and the boundary map is the composition*

$$\mathbb{Q}[-1] \to \mathbb{Q}/\mathbb{Z}[-1] \xrightarrow{e} \mathbb{Q}/\mathbb{Z}[0] \xrightarrow{\beta} \mathbb{Z}[1]. \tag{3}$$

*Proof.* Let $\alpha_m$ be the map $\oplus_{i=0}^{m-1} \mathbb{Z} \to \mathbb{Q}/\mathbb{Z} \oplus \mathbb{Q}$ which sends $x^{(i)}$ to $(\frac{x}{2} + \frac{ix}{m} + \mathbb{Z}, \frac{x}{m})$. This is a map of $\hat{G}$-modules, because for $j \in \hat{G}/m$, we have

$$\alpha_m(j \cdot x^{(i)}) = \alpha_m(x^{(i+j)}) = (\tfrac{x}{2} + \tfrac{x(j+i)}{m} + \mathbb{Z}, \tfrac{x}{m}) = j \cdot (\tfrac{x}{2} + \tfrac{ix}{m} + \mathbb{Z}, \tfrac{x}{m}) = j \cdot \alpha_m(x^{(i)}).$$

The maps $\alpha_m$ are compatible with $\delta^n_m$ (here we need the correcting factor $\frac{x}{2}$):

$$\alpha_{nm}(\delta^n_m(x^{(i)})) = \alpha_{nm}\left(\sum_{j=0}^{n-1} x^{(i+jm)}\right) = \sum_{j=0}^{n-1} \left(\tfrac{x}{2} + \tfrac{x(i+jm)}{nm} + \mathbb{Z}, \tfrac{x}{nm}\right)$$
$$= \left(\tfrac{nx}{2} + \tfrac{xi}{m} + \tfrac{x(n-1)}{2} + \mathbb{Z}, \tfrac{x}{m}\right) = \left(\tfrac{x}{2} + \tfrac{xi}{m} + \mathbb{Z}, \tfrac{x}{m}\right) = \alpha_m(x^{(i)}).$$

Hence we get a map $\mathrm{colim}_m \oplus_{i=1}^{m-1} \mathbb{Z} \xrightarrow{\alpha} M$, which is comptible with the action of $\hat{G}$. Consider the following diagram of maps of $\hat{G}$-modules with exact rows, where $S'_m = \frac{1}{m} S_m$ and $S' = \mathrm{colim}_m S'm$.

$$\begin{array}{ccccccc}
\mathbb{Z} & \xrightarrow{\Delta} & \mathrm{colim}_m \oplus_{i=0}^{m-1} \mathbb{Z} & \xrightarrow{t-1} & \mathrm{colim}_m \oplus_{i=0}^{m-1} \mathbb{Z} & \xrightarrow{S'} & \mathbb{Q} \\
\| & & S' \downarrow & & \alpha \downarrow & & \| \\
\mathbb{Z} & \to & \mathbb{Q} & \xrightarrow{g} & M & \xrightarrow{p} & \mathbb{Q}
\end{array}$$

It is easy to see that the outer squares commute. On the other hand,

$$\alpha_m(t-1)x^{(i)} = \alpha_m(x^{(i+1)}) - \alpha_m(x^{(i)}) = (\tfrac{x}{m} + \mathbb{Z}, 0) = g(\tfrac{x}{m}) = gS'(x^{(i)}).$$

□

**Proposition 4.4.** *If $\mathcal{G}^\cdot$ is a complex of $\mathbb{Q}$-vector spaces $\mathcal{T}_{\hat{G}}$, then $R\gamma_*\gamma^*\mathcal{G}^\cdot = \mathcal{G}^\cdot \oplus \mathcal{G}^\cdot[-1]$, and the cup product map*

$$R\gamma_*\gamma^*\mathcal{G}^\cdot \cong \mathcal{G}^\cdot \oplus \mathcal{G}^\cdot[-1] \xrightarrow{e_{\gamma^*\mathcal{G}^\cdot}} \mathcal{G}^\cdot[1] \oplus \mathcal{G}^\cdot \cong R\gamma_*\gamma^*\mathcal{G}^\cdot[1]$$

*is given by multiplication by the matrix $\begin{pmatrix} 0 & 0 \\ 1 & 0 \end{pmatrix}$.*

*Proof.* By theorems 3.3 and 4.3 we have

$$(\mathbb{Q} \oplus \mathbb{Q}[-1]) \otimes \mathcal{G}^\cdot \cong R\gamma_*\mathbb{Q} \otimes \mathcal{G}^\cdot \xrightarrow[\sim]{\tau} R\gamma_*\gamma^*\mathcal{G}^\cdot.$$

To calculate the multiplication by $e$, consider the following diagram, which is commutative by naturality of $e$:

$$\begin{array}{ccc} (R\gamma_*\mathbb{Q}) \otimes \mathcal{G}^\cdot & \xrightarrow{R\gamma_* e_\mathbb{Q} \otimes \mathrm{id}} & (R\gamma_*\mathbb{Q}[1]) \otimes \mathcal{G}^\cdot \\ \tau \downarrow \sim & & \tau \downarrow \sim \\ R\gamma_*\gamma^*\mathcal{G}^\cdot & \xrightarrow{R\gamma_* e_{\gamma^*\mathcal{G}^\cdot}} & R\gamma_*\gamma^*\mathcal{G}^\cdot[1] \end{array}$$

So it suffices to calculate the map

$$R\gamma_*\mathbb{Q} \cong \mathbb{Q} \oplus \mathbb{Q}[-1] \xrightarrow{e} \mathbb{Q}[1] \oplus \mathbb{Q} \cong R\gamma_*\mathbb{Q}[1].$$

i.e. calculate the boundary map in the short exact sequence $0 \to R\gamma_*\mathbb{Q} \to R\gamma_*\mathbb{Q} \otimes N \to R\gamma_*\mathbb{Q} \to 0$. It clearly suffices to do this in the copy indexed my $m = 1$. In this case, it follows from the diagram

$$\begin{array}{ccccccccc} 0 & \longrightarrow & \mathbb{Q} & \longrightarrow & N \otimes \mathbb{Q} & \longrightarrow & \mathbb{Q} & \longrightarrow & 0 \\ & & t-1 \downarrow =0 & & t-1 \downarrow & & t-1 \downarrow =0 & & \\ 0 & \longrightarrow & \mathbb{Q} & \longrightarrow & N \otimes \mathbb{Q} & \longrightarrow & \mathbb{Q} & \longrightarrow & 0 \end{array}$$

A element $r \in \mathbb{Q}$ in the upper right corner corresponds to $(0, r) \in N \otimes \mathbb{Q}$ which maps to $(r, 0)$ under $t - 1$, and this corresponds to the element $r \in \mathbb{Q}$. □

## 5. Weil motivic cohomology

Consider the following commutative diagram of functors

$$\begin{array}{ccc} \mathcal{T}_G & \xrightarrow{\gamma_*} & \mathcal{T}_{\hat{G}} \\ \Gamma_{\bar{X}} \downarrow & & \Gamma_{\bar{X}} \downarrow \\ \mathrm{Mod}_G & \xrightarrow{\gamma_*} & \mathrm{Mod}_{\hat{G}} \\ \Gamma_G \downarrow & & \Gamma_{\hat{G}} \downarrow \\ \mathrm{Ab} & = & \mathrm{Ab} \end{array}$$

For a sheaf $\mathcal{F}$ in $\mathcal{T}_G$, Lichtenbaum defines Weil-étale motivic cohomology $H^i_W(X, \mathcal{F})$ as the derived functors of the global section functor $\Gamma_G \circ \Gamma_{\bar{X}}$. Given a sheaf $\mathcal{G}$ in $\mathcal{T}_{\hat{G}}$, we also write $H^i_W(X, \mathcal{G})$ for $H^i_W(X, \gamma^*\mathcal{G})$. Similarly, for a sheaf $\mathcal{G}$ in $\mathcal{T}_{\hat{G}}$, étale cohomology $H^i_{\text{ét}}(X, \mathcal{G})$ is the derived functor of $\Gamma_{\hat{G}} \circ \Gamma_{\bar{X}}$. For a sheaf $\mathcal{F}$ in $\mathcal{T}_G$, the étale cohomology on $\bar{X}$ is not continuous in general, and we define the $\hat{G}$-continuous cohomology $H^i_{\hat{G}}(\bar{X},\mathcal{F})$ to be the derived functors of $\gamma_* \circ \Gamma_{\bar{X}}$. We get the following spectral sequences for composition of functors

$$H^s(G, H^t_{\text{ét}}(\bar{X},\mathcal{F})) \Rightarrow H^{s+t}_W(X,\mathcal{F}) \qquad (4)$$

$$H^s_{cont}(\hat{G}, H^t_{\hat{G}}(\bar{X},\mathcal{F})) \Rightarrow H^{s+t}_W(X,\mathcal{F}) \qquad (5)$$

$$H^s_{cont}(\hat{G}, H^t_{\text{ét}}(\bar{X},\mathcal{G})) \Rightarrow H^{s+t}_{\text{ét}}(X,\mathcal{G}). \qquad (6)$$

Since $\mathbb{Z}$ has cohomological dimension 1, the former spectral sequence breaks up into short exact sequence

$$0 \to H^{t-1}_{\text{ét}}(\bar{X},\mathcal{F})_G \to H^t_W(X,\mathcal{F}) \to H^t_{\text{ét}}(\bar{X},\mathcal{F})^G \to 0. \qquad (7)$$

Because the Hochschild-Serre spectral sequence is multiplicative, we get the following commutative diagram

$$\begin{array}{ccc} H^t_W(X,\mathcal{F}) & \longrightarrow & H^t(\bar{X},\mathcal{F})^G \\ e \downarrow & & can \downarrow \\ H^{t+1}_W(X,\mathcal{F}) & \longleftarrow & H^t(\bar{X},\mathcal{F})_G \end{array}$$

The identity $\Gamma_G \circ \Gamma_{\bar{X}} \circ \gamma_* = \Gamma_G \circ \Gamma_{\bar{X}}$ gives another spectral sequence

$$H^s_{\text{ét}}(X, R^t\gamma_*\mathcal{F}) \Rightarrow H^{s+t}_W(X,\mathcal{F}).$$

5.1. **Continuous Weil cohomology.** For a pro-system of étale sheaves, Jannsen [9] defined continuous cohomology as the derived functors of $\lim \circ \Gamma_X$. There is an analog for the Weil-étale topology.

**Lemma 5.1.** *Inverse limits in the categories $\mathcal{T}_G$ and $\mathcal{T}_{\hat{G}}$ exist. More precisely, $\lim_G$ agrees with the inverse limit in the category of sheaves of abelian groups, whereas $\lim_{\hat{G}} = \gamma_* \lim \gamma^*$. In particular, $R\lim_{\hat{G}} = R\gamma_* R\lim \gamma^*$.*

*Proof.* For an inverse system $\mathcal{G}^{\cdot}$ of sheaves in $\mathcal{T}_{\hat{G}}$, and a sheaf $\mathcal{G}'$ in $\mathcal{T}_{\hat{G}}$ we get by adjointness and $\mathcal{G}^i = \gamma_*\gamma^*\mathcal{G}^i$

$$\text{Hom}_{\hat{G}}(\mathcal{G}', \gamma_* \lim \gamma^*\mathcal{G}^{\cdot}) = \text{Hom}_G(\gamma^*\mathcal{G}', \lim \gamma^*\mathcal{G}^{\cdot})$$
$$= \lim \text{Hom}_G(\gamma^*\mathcal{G}', \gamma^*\mathcal{G}^{\cdot}) = \lim \text{Hom}_{\hat{G}}(\mathcal{G}', \mathcal{G}^{\cdot}) = \text{Hom}_{\hat{G}}(\mathcal{G}', \lim_{\hat{G}} \mathcal{G}^{\cdot}).$$
□

By lemma 5.1 there is a commutative diagram of topoi

$$\begin{array}{ccc} \mathcal{T}^{\mathbb{N}}_G & \xrightarrow{\gamma_*} & \mathcal{T}^{\mathbb{N}}_{\hat{G}} \\ \lim \downarrow & & \lim_{\hat{G}} \downarrow \\ \mathcal{T}_G & \xrightarrow{\gamma_*} & \mathcal{T}_{\hat{G}} \end{array}$$

The functor $\lim_{\hat{G}}$ on $\mathcal{T}^{\mathbb{N}}_{\hat{G}}$ corresponds to the functor $\lim$ on $(Shv/X)_{\text{ét}}$ under the identification of sheaves on $X$ with sheaves on $\bar{X}$ with a $\hat{G}$-action. Let $\mathcal{F}^{\cdot} \in \mathcal{T}^{\mathbb{N}}_G$ be a pro-system of sheaves of $G$-modules on $\bar{X}$. We define $H^i(\bar{X},(\mathcal{F}^{\cdot}))$ as the derived functors of $\lim_G \circ \Gamma_{\bar{X}}$. This is usually called continuous cohomology, but note that these cohomology groups are not continuous $\hat{G}$-modules in general. Because of lemma 5.1, $H^i(\bar{X},(\mathcal{F}^{\cdot}))$ can be calculated by disregarding the $G$-action. On the other hand, we define $H^i(X,(\mathcal{F}^{\cdot}))$ to be the derived functors of $\Gamma_G \circ \lim_G \circ \Gamma_{\bar{X}}$. If $(\mathcal{F}^{\cdot}) = (\gamma^*\mathcal{G}^{\cdot})$ for $(\mathcal{G}^{\cdot}) \in \mathcal{T}^{\mathbb{N}}_{\hat{G}}$, then this agrees with continuous cohomology in the sense of Jannsen, because $\Gamma_G \circ \lim_G \circ \gamma^* = \Gamma_{\hat{G}} \circ \gamma_* \circ \lim_G \circ \gamma^* = \Gamma_{\hat{G}} \circ \lim_{\hat{G}}$. Finally, we define $\hat{G}$-continuous cohomology $H^*_{\hat{G}}(\bar{X}, (\mathcal{F}^{\cdot}))$ as the derived functors of $\gamma_* \circ \lim_G \circ \Gamma_{\bar{X}}$. It is a continuous $\hat{G}$-module for the *discrete* topology. There is a map of spectral sequences

$$\begin{array}{ccc} H^s_{cont}(\hat{G}, H^t_{\hat{G}}(\bar{X}, (\mathcal{F}^{\cdot}))) & \Rightarrow & H^{s+t}(X, (\mathcal{F}^{\cdot})) \\ \downarrow & & \parallel \\ H^s(G, H^t(\bar{X}, (\mathcal{F}^{\cdot}))) & \Rightarrow & H^{s+t}(X, (\mathcal{F}^{\cdot})) \end{array} \qquad (8)$$

Note that this holds without any assumptions on the cohomology groups. The upper spectral sequence differs from Jannsen's spectral sequence [9, cor. 3.4] for finitely generated cohomology groups

$$H^s_{cont}(\hat{G}, H^t(\bar{X}, (\mathcal{F}^{\cdot}))) \Rightarrow H^{s+t}(X, (\mathcal{F}^{\cdot})),$$

because Jannsen considers the coefficients with the limit topology, whereas we always work with the discrete topology.

**Lemma 5.2.** *a) If $(\mathcal{F}^{\cdot})$ is an inverse system of sheaves in $\mathcal{T}_G$, then the sheaf $R^s \lim \mathcal{F}^{\cdot}$ is the sheaf to the presheaf which sends an étale $U$ over $\bar{X}$ to $H^s(U, (\mathcal{F}^{\cdot}))$.*
*b) If $(\mathcal{G}^{\cdot})$ is an inverse system of sheaves in $\mathcal{T}_{\hat{G}}$, then the sheaf $R^s \lim_{\hat{G}} \mathcal{G}^{\cdot}$ is the sheaf to the presheaf $U \mapsto H^s_{\hat{G}}(U, (\gamma^* \mathcal{G}^{\cdot}))$.*

*Proof.* a) follows with the same proof as b) by erasing all $\gamma_*$ and $\gamma^*$.

b) Let $\mathcal{G}^{\cdot} \to I^{\cdot}_*$ be an injective resolution of inverse systems. Then $R^s \lim_{\hat{G}}(\mathcal{G}^{\cdot})$ is by definition $\mathcal{H}^s(\gamma_* \lim \gamma^* I^{\cdot}_*) = a\mathcal{H}^s(i\gamma_* \lim \gamma^* I^{\cdot}_*)$, where $a$ is the sheafification functor and $i$ the inclusion of sheaves into presheaves. On the other hand, for every pro-system of sheaves $(\lim \mathcal{G}^{\cdot})(U) = \lim(\mathcal{G}^{\cdot}(U))$, for every sheaf $(\gamma_* \mathcal{G})(U) = \gamma_*(\mathcal{G}(U))$ and $\gamma^* \mathcal{G}(U) = \mathcal{G}(U)$, hence

$$\mathcal{H}^s(i\gamma_* \lim \gamma^* I^{\cdot}_*)(U) = H^s(\gamma_* \lim I^{\cdot}_*(U)) =: H^s_{\hat{G}}(U, (\gamma^* \mathcal{G}^{\cdot})).$$

□

We will be especially interested in the inverse system $\mathbb{Z}/l^r(n)$. Since $R\gamma_*(\mathbb{Z}/l^r(n)) = (\mathbb{Z}/l^r(n))$, we get

$$R\gamma_* R \lim \mathbb{Z}/l^r(n) \cong R \lim_{\hat{G}} \mathbb{Z}/l^r(n).$$

To be consistent with the literature, we write $H^i_{cont}(X, \mathbb{Z}_l(n))$ for $H^i(X, (\mathbb{Z}/l^r(n)))$.

## 6. Motivic cohomology

From now on we will assume that $X$ is a smooth scheme over $\mathbb{F}_q$. For $n \geq 0$, let $z^n(X, *)$ be the cycle complex of Bloch [1], and let $\mathbb{Z}(n)$ be the complex of étale sheaves sending an étale $U$ over $X$ to $z^n(U, 2n - *)$. More generally, given an abelian group $A$ we define the complex of étale sheaves $A(n)$ to be $\mathbb{Z}(n) \otimes A$, and $H^i_{\text{ét}}(X, A(n))$ as the hypercohomology of this complex. In order to make our formulas work in general, we define

$$\mathbb{Z}(n) = \mathbb{Q}/\mathbb{Z}'(n)[-1] \quad \text{for } n < 0,$$

where $\mathbb{Q}/\mathbb{Z}'$ is the prime to $p$-part of $\mathbb{Q}/\mathbb{Z}$ (although this is not "motivic cohomology" in the sense of extensions of motives). Then for any $n$, there is are quasi-isomorphisms of étale sheaves [6, 7], where $\nu^n_{p^r} = W_r \Omega^n_{X, \log}$ is the logarithmic de

Rham-Witt sheaf

$$\mathbb{Z}/m(n) \cong \begin{cases} \mu_m^{\otimes n}[0] & p \nmid m \\ \nu_{p^r}^n[-n] & m = p^r \end{cases}.$$

Furthermore, $H^i_{\text{ét}}(U, \mathbb{Q}(n)) = H^i(U, \mathbb{Q}(n))$, where the right hand side is motivic cohomology, i.e. the hypercohomology of $\mathbb{Q}(n)$ for the Zariski cohomology. Theorem 3.4 shows that there is a quasi-isomorphism of complexes

$$R\gamma_*\mathbb{Z} \otimes^L \mathbb{Z}(n) \xrightarrow{\sim} R\gamma_*\mathbb{Z}(n).$$

In particular, we get the following consequences

**Theorem 6.1.** *There is a distinguished triangle in $\mathcal{T}_{\hat{G}}$,*

$$\mathbb{Z}(n) \to R\gamma_*\gamma^*\mathbb{Z}(n) \to \mathbb{Q}(n)[-1] \xrightarrow{\delta} \mathbb{Z}(n)[1].$$

*In particular, there is a long exact sequence*

$$\ldots \to H^i_{\text{ét}}(X, \mathbb{Z}(n)) \to H^i_W(X, \mathbb{Z}(n)) \to H^{i-1}_{\text{ét}}(X, \mathbb{Q}(n)) \xrightarrow{\delta} H^{i+1}_{\text{ét}}(X, \mathbb{Z}(n)) \to \ldots$$

*where the map $\delta$ is the composition*

$$H^{i-1}_{\text{ét}}(X, \mathbb{Q}(n)) \to H^{i-1}_{\text{ét}}(X, \mathbb{Q}/\mathbb{Z}(n)) \xrightarrow{\cdot e} H^i_{\text{ét}}(X, \mathbb{Q}/\mathbb{Z}(n)) \xrightarrow{\beta} H^{i+1}_{\text{ét}}(X, \mathbb{Z}(n)).$$

**Theorem 6.2.** *Let $X$ be a smooth variety over $\mathbb{F}_q$, then*

$$H^i_W(X, \mathbb{Z}/m(n)) \cong H^i_{\text{ét}}(X, \mathbb{Z}/m(n))$$
$$H^i_W(X, \mathbb{Q}(n)) \cong H^i(X, \mathbb{Q}(n)) \oplus H^{i-1}(X, \mathbb{Q}(n)),$$

*The cup product with $e \in \text{Ext}^1_G(\mathbb{Z}, \mathbb{Z})$ is multiplication with the matrix $\begin{pmatrix} 0 & 0 \\ 1 & 0 \end{pmatrix}$.*

*Proof.* The statement of torsion sheaves is corollary 3.4, and rational statement is proposition 4.4. □

**Lemma 6.3.** *Let $\bar{X}$ be a smooth variety of dimension $d$ over an algebraically closed field of characteristic $p$, and let $l \neq p$.*

*a) If $n \geq d$, then $H^i_{\text{ét}}(\bar{X}, \mathbb{Z}(n)) \cong H^i(\bar{X}, \mathbb{Z}(n))$, and the latter group is zero for $i > n + d$.*

*b) If $n < d$, then $H^i_{\text{ét}}(\bar{X}, \mathbb{Z}(n))$ is zero for $i > 2d+1$, torsion for $i > 2n$, $p$-torsion free for $i > n + d + 1$, $p$-divisible for $i > n + d$ and $l$-divisible for $i > 2d$.*

*Proof.* a) If $\epsilon_* : \bar{X}_{\text{ét}} \to \bar{X}_{Zar}$ is the change of topology map, then $R\epsilon_*\mathbb{Q}(n)_{\text{ét}} \cong \mathbb{Q}(n)$ for any $n$, see [5]. For $n \geq d$ and $p \nmid m$ it follows from Suslin [20] that $R\epsilon_*\mathbb{Z}/m(n)_{\text{ét}} \cong \mathbb{Z}/m(n)$. On the other hand, by [6], $R\epsilon_*\mathbb{Z}/p(n)_{\text{ét}} \cong \mathbb{Z}/p(n)$, because both sides are zero for $n > d$, and because $R^1\gamma_*\nu^d = 0$. Hence we have $H^i_{\text{ét}}(\bar{X}, \mathbb{Z}(n)) \cong H^i(\bar{X}, \mathbb{Z}(n))$ for any $i$ and $n \geq d$.

b) By comparing to Zariski cohomology we get $H^i_{\text{ét}}(X, \mathbb{Q}(n)) = 0$ for $i > 2n$. For mod $p$ coefficients, $\mathbb{Z}/p(n)_{\text{ét}} \cong \nu^n[-n]$ implies that $H^i_{\text{ét}}(\bar{X}, \mathbb{Z}/p^r(n)) = H^{i-n}_{\text{ét}}(\bar{X}, \nu^n) = 0$ for $i > d + n$, because $\text{cd}_p \bar{X} = d$. For mod $l$ coefficients, $\mathbb{Z}/l(n)_{\text{ét}} = \mu_l^{\otimes n}$ implies that $H^i_{\text{ét}}(\bar{X}, \mathbb{Z}/l(n)) = 0$ for $i > 2d$ because $\text{cd}_l \bar{X} = 2d$. The statement now follows using the short exact sequence

$$0 \to H^i(\bar{X}, \mathbb{Z}(n))/l \to H^i(\bar{X}, \mathbb{Z}/l(n)) \to {}_lH^{i+1}(\bar{X}, \mathbb{Z}(n)) \to 0.$$

□

**Theorem 6.4.** *Let $X$ be a smooth variety over a finite field of dimension $d$. Then*

$$H^i_W(X, \mathbb{Z}(n)) = 0 \quad \text{for} \quad \begin{cases} i > 2d+1 & n \le d \\ i > n+d+1 & n \ge d \end{cases}.$$

*Proof.* In view of (7) and the previous lemma, the result is clear for $n \ge d$. Similarly, for $n < d$ and $i > 2d+2$, the group in question vanishes, because $H^i_{\text{ét}}(\bar{X}, \mathbb{Z}(n)) = 0$ for $i > 2d+1$. It remains to show $H^{2d+2}_W(X, \mathbb{Z}(n)) \cong H^{2d+1}_{\text{ét}}(\bar{X}, \mathbb{Z}(n))_{\mathbb{Z}} = 0$. By the lemma, $H^{2d+1}_{\text{ét}}(\bar{X}, \mathbb{Z}(n))$ is a $p$-torsion free divisible torsion group, and this property is inherited by the quotient group $H^{2d+2}_W(X, \mathbb{Z}(n))$. On the other hand, by corollary 3.4, the limit of the surjections $H^{2d+1}_W(X, \mathbb{Z}/l^r(n)) \to {}_{l^r}H^{2d+1}_W(X, \mathbb{Z}(n))$ gives a surjection $H^{2d+1}_{\text{cont}}(X, \mathbb{Z}_l(n)) \to T_l H^{2d+2}_W(X, \mathbb{Z}(n))$. By Deligne's proof of the Weil conjectures $H^i_{\text{cont}}(X, \mathbb{Q}_l(n)) = 0$ for $i > n+d+1$, hence the former group is torsion and $T_l H^{2d+2}_W(X, \mathbb{Z}(n)) = 0$. The theorem follows with the following lemma. □

**Lemma 6.5.** *Let $A$ be a divisible $l$-torsion group such that $T_l A = 0$. Then $A = 0$.*

*Proof.* The inverse limit of the short exact sequences

$$0 \to {}_{l^r}A \to A \xrightarrow{l^r} A \to 0$$

gives a short exact sequence

$$0 \to T_l A \to \lim(A, l) \to A \to 0,$$

because by divisibility the inverse system $({}_{l^r}A)_r$ is Mittag-Leffler. If $T_l A = 0$, then $A = \lim(A, l)$ is uniquely divisible and torsion, hence zero. □

**Theorem 6.6.** *Let $X$ be a connected smooth projective variety of dimension $d$ over $\mathbb{F}_q$. Then there is a commutative diagram with indicated isomorphisms*

$$\begin{array}{ccccc}
H^{2d}_W(X, \mathbb{Z}(d)) & \longrightarrow & H^{2d}_{\text{ét}}(\bar{X}, \mathbb{Z}(d))^{\mathbb{Z}} & \xrightarrow{(deg,?)} & \mathbb{Z} \oplus (finite) \\
{\scriptstyle e}\downarrow & & \downarrow & & \downarrow \\
H^{2d+1}_W(X, \mathbb{Z}(d)) & \xleftarrow{\sim} & H^{2d}_{\text{ét}}(\bar{X}, \mathbb{Z}(d))_{\mathbb{Z}} & \xrightarrow{deg} & \mathbb{Z}
\end{array}$$

*Proof.* By the lemma $H^{2d}_{\text{ét}}(\bar{X}, \mathbb{Z}(d)) \cong \text{CH}^d(\bar{X})$. The kernel of the degree map $\text{CH}^d(\bar{X}) \xrightarrow{deg} \mathbb{Z}$ is divisible, and by Roitman's theorem it agrees with the $\bar{\mathbb{F}}_q$-torsion points $A(\bar{\mathbb{F}}_q)$ of the Albanese variety $A$ of $X$. Since $A(\bar{\mathbb{F}}_q)^{\mathbb{Z}} = A(\mathbb{F}_q)$ is finite, we get $A(\bar{\mathbb{F}}_q)_{\mathbb{Z}} = 0$. The statement now follows from the short exact sequence (7). □

## 7. Comparison to $l$-adic cohomology

Fix a smooth projective variety $X$ over $\mathbb{F}_q$ and an integer $n$. There are two fundamental conjectures on Weil motivic cohomology. The first one is due to Lichtenbaum:

**Conjecture 7.1. L(X,n)** *For any $i$, the group $H^i_W(X, \mathbb{Z}(n))$ is finitely generated.*

The second fundamental conjecture relates Weil cohomology to $l$-adic cohomology. The compatible homomorphisms $\mathbb{Z}(n) \otimes \mathbb{Z}_l \to \mathbb{Z}/l^r(n)$ in the derived category of sheaves of $\mathcal{T}_G$ induce a morphism $\mathbb{Z}(n) \otimes \mathbb{Z}_l \to R \lim \mathbb{Z}/l^r(n)$, hence in view of lemma 5.1 upon applying $R\gamma_*$ a map $R\gamma_* \mathbb{Z}(n) \otimes \mathbb{Z}_l \to R \lim \mathbb{Z}/l^r(n)$ in the derived category of sheaves of $\mathcal{T}_{\hat{G}}$.

**Conjecture 7.2. C(X,n)** *For every prime $l$, and any $i$, there is an isomorphism*
$$H^i_W(X, \mathbb{Z}(n)) \otimes \mathbb{Z}_l \xrightarrow{\sim} H^i_{cont}(X, \mathbb{Z}_l(n))$$

In particular, the latter conjecture implies that Weil motivic cohomology gives an integral model for $l$-adic and $p$-adic cohomology. Finally, there is the classical conjecture, due to Tate (part 1,2 for $l \neq p$) and Beilinson (part 3):

**Conjecture 7.3. T(X,n)** *For any prime number $l$,*
1. *The cycle map $CH^n(X) \otimes \mathbb{Q}_l \to H^{2n}_{cont}(X, \mathbb{Q}_l(n))$ is surjective.*
2. *The $\hat{G}$-module $H^{2n}_{cont}(\bar{X}, \mathbb{Q}_l(n))$ is semi-simple at $1$.*
3. *Rational and numerical equivalence agree with rational coefficients.*

The following theorem uses ideas from Kahn [10, 12], Soulé [19] and the author [3].

**Theorem 7.4.** *Let $X$ be a smooth projective variety over $\mathbb{F}_q$, and $n$ an integer. Then*
$$C(X, n) + C(X, d-n) \Rightarrow L(X, n) \Rightarrow C(X, n) \Rightarrow T(X, n).$$

*Proof.* To show that conjecture $C(X, n)$ and $C(X, d-n)$ imply $L(X, n)$, we repeat the argument of Kahn [12]. Since by Gabber [2] the groups $H^i_{cont}(X, \mathbb{Z}_l(n))$ have torsion only for finitely many $l$, the same holds for $H^i_W(X, \mathbb{Z}(n))$, and it suffices to show that $\bar{H}^i_W(X, \mathbb{Z}(n)) := H^i_W(X, \mathbb{Z}(n))/tor$ is finitely generated. By hypothesis and proposition 6.6, the duality pairing
$$\bar{H}^i_W(X, \mathbb{Z}(n)) \times \bar{H}^{2d+1-i}_W(X, \mathbb{Z}(d-n)) \to \bar{H}^{2d+1}_W(X, \mathbb{Z}(d)) \cong \mathbb{Z}$$
becomes a perfect pairing after tensoring with $\mathbb{Z}_l$ for any $l$, and both terms are finitely dimensional vector spaces after tensoring with $\mathbb{Q}$. In this situation, [12, lemma 3.8] shows that $\bar{H}^i_W(X, \mathbb{Z}(n))$ is finitely generated.

To show that $L(X, n)$ implies $C(X, n)$, we consider the following commutative diagram
$$\begin{array}{ccccc}
H^i_W(X, \mathbb{Z}(n)) \otimes \mathbb{Z}_l/l^r & \longrightarrow & H^i_W(X, \mathbb{Z}/l^r(n)) & \longrightarrow & {}_{l^r}H^{i+1}_W(X, \mathbb{Z}(n)) \otimes \mathbb{Z}_l \\
\downarrow & & \parallel & & \downarrow \\
H^i_{cont}(X, \mathbb{Z}_l(n))/l^r & \longrightarrow & H^i_{\text{ét}}(X, \mathbb{Z}/l^r(n)) & \longrightarrow & {}_{l^r}H^{i+1}_{cont}(X, \mathbb{Z}_l(n))
\end{array}$$
We proceed by induction on $i$. If $H^{i-1}_W(X, \mathbb{Z}(n)) \otimes \mathbb{Z}_l = H^{i-1}_{cont}(X, \mathbb{Z}_l(n))$, then the torsion subgroups of $H^i_W(X, \mathbb{Z}(n)) \otimes \mathbb{Z}_l$ and of $H^i_{cont}(X, \mathbb{Z}_l(n))$ agree. On the other hand, if these finitely generated groups had different rank, then by the snake lemma the torsion subgroup of $H^{i+1}_W(X, \mathbb{Z}(n)) \otimes \mathbb{Z}_l$ would be infinite, a contradiction.

To show that conjecture $C(X, n)$ implies $T(X, n)$, consider the diagram
$$\begin{array}{ccccc}
H^{2n}_{\text{ét}}(X, \mathbb{Z}(n)) \otimes \mathbb{Q}_l & \xrightarrow{i} & H^{2n}_W(X, \mathbb{Z}(n)) \otimes \mathbb{Q}_l & \xrightarrow{\sim} & H^{2n}_{cont}(X, \mathbb{Q}_l(n)) \\
\parallel & & e \downarrow & & e \downarrow \\
H^{2n}_{\text{ét}}(X, \mathbb{Z}(n)) \otimes \mathbb{Q}_l & \xleftarrow{\sim} & H^{2n+1}_W(X, \mathbb{Z}(n)) \otimes \mathbb{Q}_l & \xrightarrow{\sim} & H^{2n+1}_{cont}(X, \mathbb{Q}_l(n))
\end{array} \quad (9)$$

The left square commutes by theorem 6.2, and the isomorphisms hold by hypothesis, theorem 6.2 and because $H^{2n+1}_{\text{ét}}(X,\mathbb{Z}(n)) \otimes \mathbb{Q}_l = 0$. Since all groups are $\mathbb{Q}_l$-vector spaces of the same dimension, the injection $i$ is an isomorphism, hence multiplication by $e$ is an isomorphism, which implies semi-simplicity of $H^{2n}_{cont}(X,\mathbb{Q}_l(n))$. This, together with the surjectivity of the upper row, implies the strong form of Tate's conjecture [22, theorem 2.9]. Finally, the injectivity of the upper row shows that rational and numerical equivalence agree. □

*Remark.* Conjecture $C(X,n)$ for all smooth and proper varieties implies the same statement for all smooth varieties as long as $l \neq p$. This follows by using localization sequences for $l$-adic cohomology and Weil-cohomology, and de Jong's theorem on alterations, see [4, lemma 4.1] or [10, section 5] for details. On the other hand we have the following

*Example.* Let $\mathbb{A} = \mathbb{A}^1_{\mathbb{F}_q}$ be the affine line. Then Artin Schreier theory shows that $H^1_{\text{ét}}(\mathbb{A},\mathbb{Z}/p)$ is an infinite direct sum of copies of $\mathbb{Z}/p$. An easy diagram chase shows that $H^1_{\text{ét}}(\mathbb{A},\mathbb{Q}/\mathbb{Z}) \cong H^2_{\text{ét}}(\mathbb{A},\mathbb{Z})$ contains infinitely many copies of $\mathbb{Q}_p/\mathbb{Z}_p$, which gives an infinite number of copies of $\mathbb{Q}_p/\mathbb{Z}_p$ in $H^2_W(\mathbb{A},\mathbb{Z})$. This is why Lichtenbaum considers Weil cohomology with compact support in [14].

We have the following converse of theorem 7.4.

**Theorem 7.5.** *If $T(X,n)$ holds for all smooth and projective varieties $X$ over $\mathbb{F}_q$ and all $n$, then $C(X,n)$ and $L(X,n)$ hold for all $X$ and $n$.*

*Proof.* By theorem 6.2 it suffices to show the result rationally. Let $X$ be a smooth projective variety. By [3] the hypothesis implies $H^i(X,\mathbb{Q}(n)) \cong H^i_{\text{ét}}(X,\mathbb{Q}(n)) = 0$ for $i \neq 2n$, hence $H^i_W(X,\mathbb{Z}(n)) \otimes \mathbb{Q} = 0$ for $i < 2n$ and for $i > 2n+1$. By hypothesis, the composition

$$H^{2n}_{\text{ét}}(X,\mathbb{Z}(n)) \otimes \mathbb{Q}_l \to H^{2n}_W(X,\mathbb{Z}(n)) \otimes \mathbb{Q}_l \to H^{2n}_{cont}(X,\mathbb{Q}_l(n))$$

is surjective. Since $H^{2n-1}_{\text{ét}}(X,\mathbb{Q}(n)) = 0$ the first map is an isomorphism, and by equality of rational and numerical equivalence the composition is injective. The statement in degree $2n+1$ follows from semi-simplicity, using diagram (9). □

By Tate [21, theorem 2.9] it suffices to know the statement of $T(X,n)$ for only one $l$, hence the same holds for $C(X,n)$ and $L(X,n)$. See Kahn [10] for a large number of consequences of conjecture 7.2.

*Remark.* For $l \neq p$, a more abstract formulation of (2), equivalent to the main conjecture in Kahn [10], is that there is a quasi-isomorphism in the derived category of the category $(Sm/\mathbb{F}_q)_{\text{ét}}$ of smooth schemes over $\mathbb{F}_q$ equipped with the étale topology

$$R\gamma_*\mathbb{Z}(n) \otimes \mathbb{Z}_l \xrightarrow{\sim} R\lim_{\hat{G}} \mathbb{Z}/l^r(n).$$

Here $\mathbb{Z}(n)$ is Kahn's functorial version of the cycle complex [11]. However, the complexes $\mathbb{Z}(n) \otimes \mathbb{Z}_l$ and $R\lim \mathbb{Z}/l^r(n)$ are not quasi-isomorphic as objects in the derived category of $(Sm/\mathbb{F}_q)_{W-\text{ét}}$, the category of smooth schemes over $\mathbb{F}_q$ equipped with the Weil-étale topology. For example, $\mathcal{H}^2(\mathbb{Z}(1) \otimes \mathbb{Z}_l) = 0$, and we claim that $R^2\lim \mathbb{Z}/l^r(1) \neq 0$. In view of lemma 5.2, it suffices to show that some stalk of the presheaf $U \mapsto H^2_{cont}(U,\mathbb{Z}_l(1))$ is non-trivial. Since $H^2_{cont}(U,\mathbb{Z}_l(1))$ surjects onto $T_l \operatorname{Br} U$, it suffices to show that some stalk of the presheaf $U \mapsto T_l \operatorname{Br} U$ is non-trivial. Let $X$ be a smooth scheme over $\mathbb{F}_q$ such that $\operatorname{Br} \bar{X}$ is non-trivial. If $K$ is the function field of $X$, then $\operatorname{Br} X$ injects into $\operatorname{Br} K$, and since $T_l$ is left exact, $T_l \operatorname{Br} X$ injects into $T_l \operatorname{Br} K$. If $K'/K$ is a finite Galois extension of degree $e$, then the composition $T_l \operatorname{Br} K \to T_l \operatorname{Br} K' \to T_l \operatorname{Br} K$ is multiplication by $e$, hence injective since Tate modules are torsion-free. Hence the stalk of $T_l \operatorname{Br}(-)$ at the generic point of $X$ is non-trivial.

In view of the spectral sequences (5) and (8), conjecture $C(X, n)$ follows from isomorphisms $H^i_{\hat{G}}(\bar{X}, \mathbb{Z}(n)) \otimes \mathbb{Z}_l \cong H^i_{\hat{G}}(\bar{X}, (\mathbb{Z}/l^{\cdot}(n)))$ for any $i, l$ and any smooth projective variety $\bar{X}$ and for $i = 2n, 2n+1$. These isomorphisms fit into a commutative diagram of short exact sequences of $\hat{G}$-modules

$$\begin{array}{ccccc}
R^1\gamma_* H^{i-1}_{\text{ét}}(\bar{X}, \mathbb{Z}(n)) \otimes \mathbb{Z}_l & \longrightarrow & H^i_{\hat{G}}(\bar{X}, \mathbb{Z}(n)) \otimes \mathbb{Z}_l & \longrightarrow & \gamma_* H^i_{\text{ét}}(\bar{X}, \mathbb{Z}(n)) \otimes \mathbb{Z}_l \\
\downarrow & & \downarrow & & \downarrow \\
R^1\gamma_* H^{i-1}_{cont}(\bar{X}, \mathbb{Z}_l(n)) & \longrightarrow & H^i_{\hat{G}}(\bar{X}, (\mathbb{Z}/l^{\cdot}(n))) & \longrightarrow & \gamma_* H^i_{cont}(\bar{X}, \mathbb{Z}_l(n))
\end{array}$$

Thus $C(X, n)$ follows from isomorphisms of $\hat{G}$-modules

$$H^i_{\text{ét}}(\bar{X}, \mathbb{Z}(n)) \otimes \mathbb{Z}_l = \gamma_* H^i_{\text{ét}}(\bar{X}, \mathbb{Z}(n)) \otimes \mathbb{Z}_l \xrightarrow{\sim} \gamma_* H^i_{cont}(\bar{X}, \mathbb{Z}_l(n))$$

$$H^{i-1}_{\text{ét}}(\bar{X}, \mathbb{Z}(n)) \otimes \mathbb{Q}_l = R^1\gamma_* H^{i-1}_{\text{ét}}(\bar{X}, \mathbb{Z}(n)) \otimes \mathbb{Z}_l \xrightarrow{\sim} R^1\gamma_* H^i_{cont}(\bar{X}, \mathbb{Z}_l(n)).$$

It is interesting to observe that the former isomorphism was Tate's original formulation of his conjecture [23].

## 8. Values of zeta-functions, Examples

We can use the method of Milne/ Lichtenbaum to find expressions for values of zeta functions. Since $e \in \operatorname{Ext}^1_G(\mathbb{Z}, \mathbb{Z})$ satisfies $e^2 = 0$, the Weil étale cohomology groups $H^*_W(X, \mathbb{Z}(n))$ form a complex under product with $e$. For a bounded complex $C^{\cdot}$ of abelian groups with finite cohomology groups one defines

$$\chi(C^{\cdot}) := \prod_i |H^i(C^{\cdot})|^{(-1)^i}.$$

Let $X$ be a smooth projective scheme over $\mathbb{F}_q$ and $\zeta(X, s) = Z(X, q^{-s})$ be its zeta function. Following Milne [15], we let

$$\chi(X, \mathcal{O}_X, n) = \sum_{i \leq n, j \leq d} (-1)^{i+j}(n - i) \dim H^j(X, \Omega^i).$$

The conclusion of the following theorem has been proved by Lichtenbaum for $n = 0$ in [14].

**Theorem 8.1.** *Let $X$ be a smooth projective variety such that $C(X, n)$ holds. Then the order $\rho_n$ of the pole of $\zeta(X, s)$ at $s = n$ is $\operatorname{rank} H^{2n}_W(X, \mathbb{Z}(n))$, and*

$$\zeta(X, s) = \pm(1 - q^{n-s})^{-\rho_n} \cdot \chi(H^*_W(X, \mathbb{Z}(n)), e) \cdot q^{\chi(X, \mathcal{O}_X, n)} \quad \text{as } s \to n. \tag{10}$$

*If furthermore $C(X, d - n)$ holds, then the leading coefficient can also be expressed as*

$$\chi(H^*_W(X, \mathbb{Z}(n)), e) = \prod_i |H^i_W(X, \mathbb{Z}(n))_{tor}|^{(-1)^i} \cdot R^{-1} \tag{11}$$

*where $R$ is the determinant of the pairing*

$$H^{2n}_W(X, \mathbb{Z}(n)) \times H^{2(d-n)}_W(X, \mathbb{Z}(d-n)) \to H^{2d}_W(X, \mathbb{Z}(d)) \xrightarrow{deg} \mathbb{Z}.$$

*Proof.* Since the isomorphism in degrees $2n$ and $2n+1$ imply semi-simplicity of $l$-adic cohomology, the first two statements follow by comparing to the formulas for $l$-adic cohomology in [15, theorem 0.1].

To prove c), first recall that for a short exact sequence $0 \to A^{\cdot} \to B^{\cdot} \to C^{\cdot} \to 0$ we have $\chi(A^{\cdot}) \cdot \chi(C^{\cdot}) = \chi(B^{\cdot})$. Thus we have to show that the regulator term agrees with $\chi(H_W^*(X, \mathbb{Z}(n))/tor, e)$. By hypothesis, the latter complex only consist of the upper map in the following diagram

$$\begin{array}{ccc} H_W^{2n}(X, \mathbb{Z}(n))/tor & \xrightarrow{e} & H_W^{2n+1}(X, \mathbb{Z}(n))/tor \\ d \downarrow & & \| \\ \mathrm{Hom}(H_W^{2(d-n)}(X, \mathbb{Z}(d-n)), \mathbb{Z}) & \xrightarrow{e} & \mathrm{Hom}(H_W^{2(d-n)}(X, \mathbb{Z}(d-n)), H_W^{2d+1}(X, \mathbb{Z}(d))). \end{array}$$

Comparing with $l$-adic cohomology, one sees that the right vertical map is an isomorphism. The lower horizontal map is an isomorphism by theorem 6.6. Hence

$$\chi(H_W^*(X, \mathbb{Z}(n))/tor, e) = \frac{|\ker e|}{|\operatorname{coker} e|} = \frac{|\ker d|}{|\operatorname{coker} d|} = \frac{1}{R}.$$

□

8.1. **Examples.**

**Proposition 8.2.** *Conjecture $C(X, n)$ holds for $n \leq 0$. In particular, (10) holds for all $X$ and $n \leq 0$.*

*Proof.* For $n = 0$, it has been verified in Kahn [10] that $\mathbb{Z}^c \otimes \mathbb{Z}_l \cong R \lim \mathbb{Z}/l^r$. The same argument works for $l = p$. Taking étale hypercohomology of the quasi-isomorphism $R\gamma_*\mathbb{Z} \otimes \mathbb{Z}_l \cong R \lim \mathbb{Z}/l^r$ gives statement $C(X, 0)$.

For $n < 0$, the proposition follows because $H_{cont}^i(X, \mathbb{Q}_l(n)) = 0$ by Deligne's proof of the Weil conjectures, hence

$$H_W^i(X, \mathbb{Z}(n)) \otimes \mathbb{Z}_l := H_{\text{ét}}^{i-1}(X, \mathbb{Q}_l/\mathbb{Z}_l(n)) \xrightarrow[\sim]{\beta} H_{cont}^i(X, \mathbb{Z}_l(n)).$$

□

**Theorem 8.3.** *Assume that $X$ is smooth and projective and that the cycle map $\operatorname{Pic} X \otimes \mathbb{Q}_l \to H_{cont}^2(X, \mathbb{Q}_l(1))$ is surjective for some $l$. Then $C(X, 1)$ holds. In particular, (10) holds for $X$ and $n = 1$.*

*Proof.* We have $H_{\text{ét}}^i(X, \mathbb{Z}(1)) \otimes \mathbb{Q} = H^{i-1}(X, \mathbb{G}_m) \otimes \mathbb{Q} = 0$ and $H_{cont}^i(X, \mathbb{Q}_l(1)) = 0$ for $i \neq 2, 3$. Consider the diagram

$$\begin{array}{ccccc} \operatorname{Pic} X \otimes \mathbb{Q}_l & \xrightarrow{\sim} & H_W^2(X, \mathbb{Z}(1)) \otimes \mathbb{Q}_l & \xrightarrow{surj} & H_{cont}^2(X, \mathbb{Q}_l(1)) \\ \| & & e \downarrow & & e \downarrow \\ \operatorname{Pic} X \otimes \mathbb{Q}_l & \xleftarrow{\sim} & H_W^3(X, \mathbb{Z}(1)) \otimes \mathbb{Q}_l & \longrightarrow & H_{cont}^3(X, \mathbb{Q}_l(1)) \end{array}$$

In weight 1, we know that rational and numerical equivalence agree rationally, hence the upper composition is an isomorphism. On the other hand, by Milne [15, prop. 0.3], the surjectivity of the cycle map implies semi-simplicity of $H_{cont}^2(X, \mathbb{Q}_l(1))$, including $l = p$, hence the right vertical map is an isomorphism. □

In particular, the conclusion holds for Hilbert modular surfaces, Picard modular surfaces, Siegel modular threefolds, and in characteristic at least 5 for supersingular and elliptic K3 surfaces [22]. We use Soulé's method to produce more examples.

**Proposition 8.4.** *Let $X = X_1 \times \ldots \times X_d$ be a product of smooth projective curves over $\mathbb{F}_q$, and let $n \leq 1$ or $n \geq d-1$. Then $C(X,n)$ holds for $X$.*

*Proof.* In view of theorem 6.2 it suffices to show that $H^i(X, \mathbb{Z}(n)) \otimes \mathbb{Q} = 0$ for $i < 2n$, that $H^{2n}(X, \mathbb{Z}(n)) \otimes \mathbb{Q}_l \cong H^{2n}_{cont}(X, \mathbb{Q}_l(n))$, and that the Frobenius acts semi-simply at 1 on $H^{2n}_{cont}(\bar X, \mathbb{Q}_l(n))$. We essentially repeat the proof of Soulé, [19, théoréme 3], adapted to our situation.

If we write $X_i = \mathbf{1} \oplus X_i^+ \oplus \mathcal{L}$ in the category of Chow motives, then $X$ is sum of motives of the form $M = \otimes_{s=1}^j X_{n_s}^+ \otimes \mathcal{L}^k$, with $0 \leq j+k \leq d$. Such a motive $M$ has a Frobenius endomorphism $F_M$, and $F_M$ has a minimal polynomial $P_M(u)$ such that all roots of $P_M(u)$ have absolute valuation equal to $q^{\frac{j+2k}{2}}$ [19, prop. 3.1.2]. Since the Frobenius $F_M$ acts on $H^i(X, \mathbb{Z}(n)) \otimes \mathbb{Q}$ as multiplication by $q^n$ [19, prop. 1.5.2], we get $0 = P_M(F_M) = P_M(q^n)$. Hence for $j+2k \neq 2n$, multiplication by the unit $P_M(q^n)$ is zero, hence $H^i(X, \mathbb{Z}(n)) \otimes \mathbb{Q} = 0$. By our choice of $n$, $j+2k = 2n$ can only happen for $j = 0$ or $j = 2$.

If $j = 0$, then $M = \mathcal{L}^n$, and by the projective space formula for motivic cohomology we get $H^i(\mathcal{L}^n, \mathbb{Z}(n)) \otimes \mathbb{Q}_l = H^{i-2n}(\mathbb{F}_q, \mathbb{Z}(0)) \otimes \mathbb{Q}_l$. The latter group is zero except for $i = 2n$, in which case it is isomorphic to $H^{2n}_{cont}(\mathcal{L}^n, \mathbb{Q}_l(n)) = \mathbb{Q}_l$. Also, the Galois group acts trivially, in particular semi-simply on the latter group.

If $j = 2$, then $M = X^+ \otimes Y^+ \otimes \mathcal{L}^{n-1}$, and

$$H^i(M, \mathbb{Z}(n)) \otimes \mathbb{Q}_l = H^{i-2n+2}(X^+ \otimes Y^+, \mathbb{Z}(1)) \otimes \mathbb{Q}_l = H^{i-2n+1}(X^+ \otimes Y^+, \mathbb{G}_m) \otimes \mathbb{Q}_l$$

The latter group is zero for $i - 2n + 1 < 1$, i.e. for $i < 2n$, because the group of global sections of a projective variety over a finite field is finite. On the other hand,

$$H^{2n}(M, \mathbb{Z}(n)) \otimes \mathbb{Q}_l \cong \mathrm{CH}^1(X^+ \otimes Y^+) \otimes \mathbb{Q}_l$$
$$\cong H^2_{cont}(X^+ \otimes Y^+, \mathbb{Q}_l(1)) \cong H^{2n}_{cont}(M, \mathbb{Q}_l(n))$$

by Tate's theorem [21]. Furthermore, Tate's theorem also implies that the Galois group acts semi-simply at 1 on the module $H^2_{cont}(X^+ \otimes Y^+, \mathbb{Q}_l(1))$. For $l = p$, the same statement follows by [15, prop. 0.3]. □

As in Soulé, let $A(k)$ be the subclass of smooth projective varieties generated by products of curves and the following operations:

1. If $X$ and $Y$ are in $A(k)$, then $X \coprod Y$ is in $A(k)$.
2. If $Y$ is in $A(k)$, and there are morphisms $c : X \to Y$ and $c' : Y \to X$ in the category of Chow motives, such that $c' \circ c : X \to X$ is multiplication by a constant, then $X$ is in $A(k)$.
3. If $k'$ is a finite extension of $k$, and $X \times_k k'$ is in $A(k')$, then $X$ is in $A(k)$.
4. If $Y$ is a closed subscheme of $X$ and $Y$ and $X$ are in $A(k)$, then the blow-up $X'$ of $X$ along $Y$ is in $A(k)$.

**Theorem 8.5.** *Let $X$ be a variety of dimension $d$ in $A(\mathbb{F}_q)$. Then $C(X, n)$ and $L(X, n)$ hold for $n \leq 1$ or $n \geq d-1$. In particular, (10) and (11) holds for $X$ and $n \leq 1$ or $n \geq d-1$.*

*Proof.* The statement holds for products of curves, and it is clear that if $X$ and $Y$ satisfy $C(X, n)$ then $X \coprod Y$ also does. In 2) and 3), the map $H^i_W(X, \mathbb{Z}(n)) \otimes \mathbb{Z}_l \to$ $H^i_{cont}(X, \mathbb{Z}(n))$ is a direct summand for the corresponding map for $Y$ and $X \times_k k'$, respectively. Finally, if $X'$ is the blow-up of $X$ along $Y$, and $Y$ has codimension $c$ in $X$, then one has $X' = X \oplus (\oplus_{j=1}^{c-1} Y \otimes \mathcal{L}^j)$. □

In particular, the conclusion of theorem holds for abelian varieties, unirational varieties of dimension at most 3, or Fermat hypersurfaces.

UNIVERSITY OF SOUTHERN CALIFORNIA, DEPARTMENT OF MATHEMATICS, DRB, 1042 W. 36TH PLACE, LOS ANGELES, CA 90089